\documentclass[11pt]{article}
\usepackage{hyperref}
\usepackage{amssymb,amsmath, amsthm}
\usepackage{graphicx}
\usepackage{graphics}
\usepackage{psfrag}

\usepackage[all]{xy}

\theoremstyle{definition}

\newcommand{\Rb}{\mathbb{R}}

\begin{document}
 
\title{Numerical Approximation of Stochastic  Neural Fields with Delays}
\author{Pedro M. Lima (a) and Evelyn Buckwar (b)\\
(a) CEMAT - Center for Computational and Stochastic Mathematics\\
Instituto Superior T\'ecnico, Universidade de Lisboa \\
A. Rovisco Pais, 1049-001 Lisboa, Portugal\\
(b)  Institute of Stochastics\\
Johannes Kepler University\\
Altenbergerstr. 69 , 4040 Linz, Austria}
\maketitle

\begin{abstract}
We introduce  a new numerical algorithm for solving the stochastic neural field equation (NFE)
with delays. Using this algorithm we have obtained some numerical results which illustrate the
effect of noise in the dynamical behaviour of stationary solutions of the  NFE, in the presence 
of  spatially heterogeneous external inputs. 
\end{abstract}
\maketitle
%\begin{document}
\section{Introduction}
Neural Field Equations (NFE) are a powerful tool for analysing the dynamical behaviour of 
populations of neurons.  The analysis of such dynamical mechanisms is crucially important for 
understanding   a wide range of neurobiological phenomena \cite{Bressloff}.
In this work we will be concerned with the NFE in the form
\begin{equation} \label{2}
 \frac{\partial}{\partial t} V(x,t) =
I(x,t) -\alpha  V(x,t)+ \int_{\Omega} K(|x-y|) S(V(y,t-\tau(x,y)))dy,
\end{equation}
\[  t\in[0,T],  \quad  x  \in \Omega \subset \Rb,\]
where 
 $V(x,t)$ (the unknown function)  denotes the membrane potential in point $x$ at time $t$;
 $I(x,t)$ represents the external sources of excitation;
$S$ is the dependence between the firing rate of the neurons and their membrane
potentials (sigmoidal or  Heaviside function); 
and $K(|x-y|)$ gives the connectivity
between neurons at  $x$ and $y$, $\alpha$ is a constant (related to the decay rate of the potential); $\tau(x,y )>0$ is a delay, depending on the spatial variables (it results from the finite propagation speed of nervous stimulus in the brain).

 Equation \eqref{2}  (without delay) was introduced first by Wilson
and Cowan \cite{WC}, and then by Amari \cite{Amari}, to describe excitatory
and inhibitory interactions in populations of neurons.

Intensive studies of Hopf bifurcations occuring in neural fields have been carried out in the last decade.
In \cite{Bressloff} the authors investigate the occurence of spatially localized oscillations (or breathers) in two-dimensional neural fields with excitatory and inhibitory interactions;
in \cite{VF} the authors obtain sufficient conditions for the stability of stationary solutions of neural field equations; the dependency of the stationary solutions of NFE with respect to the stiffness
of the nonlinearity and the contrast of external inputs is studied in \cite{VF2}.

Though the above mentioned results were obtained analytically, numerical simulations play a fundamental role in studying brain dynamics.  Thus, the availability of efficient numerical methods is an important ingredient for improving the understanding of neural processes.
Concerning equation \eqref{2},   numerical approximations were obtained in \cite{Faye}.
The computational method applies  quadrature rule in space to reduce the problem
to a system of delay differential equations, which is then solved
by a standard algorithm for this kind of equations. A more efficient approach was recently
proposed in \cite{Hutt2010} \cite{Hutt2014}, where the authors introduce a new approach to deal with
the convolution kernel of the equation and use Fast Fourier Transforms to reduce 
significantly the computational effort required by numerical integration.
Recently, a new numerical method
for the approximation of two-dimensional neural fields has been introduced, based on an implicit 
second order scheme for the integration in time and using 
Chebyshev interpolation to reduce the dimensions of the matrices   \cite{LB}.
Some applications of this algorithm to Neuroscience problems have been discussed
in \cite{MCSI}.

As in other sciences, in Neurobiology it is well-known that better consistency with some phenomena can be provided if the effects of random processes in the system are taken into account.
In the recent work of K\"{u}hn and Riedler \cite{KR}, the authors study the effect of additive noise in Neural
Field Equations. With this purpose they introduce the stochastic integro-differential equation
\begin{equation} \label{3}
dU_t(x) =
\left(I(x,t) - \alpha  U_t(x) + \int_{\Omega} K(|x-y|) S(U_t(y))dy \right) dt + \epsilon dW_t(x),
\end{equation}
where  $t\in[0,T]$, $x \in \Omega \subset \Rb^n$, $W_t$ is a Q-Wiener process.

The main goal of the present work is to analyse the effect of noise in certain neural fields, 
which allow different types of stationary solutions. In this case we consider the following modification 
of equation \eqref{3}
\begin{equation} \label{4}
dU_t(x) =
\left(I(x,t) - \alpha U_t(x) + \int_{\Omega} K(|x-y|) S(U_{t-\tau}(y))dy \right) dt + \epsilon dW_t(x),
\end{equation}
where, as in the deterministic case,  $\tau$ is a delay, depending on the distance $|x-y|$.
Equation \eqref{4} is completed with an initial condition of the form
\begin{equation} \label{4a}
U_t(x)=U_0(x,t), \qquad t \in[-\tau_{max},0], \quad x \in \Omega,
\end{equation} 
where $U_0(x,t)$ is some given stochastic process, $\tau_{max}$  is the maximum value of the delay ( $\tau_{max}=|\Omega|/v $), where $v$ is
the propagation speed of the signals. We assume that $U_t(x)$ satisfies periodic boundary conditions in space.
We will consider domains of the form $\Omega=[-l,l]$, including the limit case when 
$ l \rightarrow \infty$.

%%%%%%%%%%%%%%%%%%%%%%%%%%%%%%
\section{Numerical Approximation}
To construct a numerical approximation of  the solution of \eqref{3} in the one-dimensional
case, we begin by expanding the solution $U_t(x)$ using the Karhunen-Loeve formula:
\begin{equation} \label{5}
U_t(x)=\sum_{k=0}^{\infty} u_t^k v_k(x),
\end{equation}
where $v_k$ are the eigenfunctions of the covariance operator of the noise in \eqref{3},
which form an orthogonal system (their explicit form is indicated below).
To derive a formula for the coefficients $u_t^k$ we take the inner product of equation \eqref{3}
with the basis functions $v_i$:
\begin{equation} \label{6}
(dU_t, v_i)=
\left[(I(x,t),v_i) - \alpha (U_t,v_i) +\left( \int_{\Omega} K(|x-y|) S(U_{t-\tau}(y))dy, v_i\right) \right] dt + \epsilon (dW_t,v_i).
\end{equation}
We expand $dW_t$ as
\begin{equation}\label{beta}
d W_t(x)=\sum_{k=0}^{\infty} v_k(x) \lambda_k  d\beta_t^k,
\end{equation}
where the functions $\beta_t^k$ form a system of independent white noises in time  and 
 $\lambda_k$ are the eigenvalues of the covariance operator of the noise.
As an important particular case, we consider the one described in \cite{KR}, p.7. 
In this case the correlation function satisfies
\[ E W_t(x) W_s(y) = \min(t,s) \frac{1}{2 \xi} \exp\left(\frac{-\pi}{4} \frac{|x-y|^2}{\xi^2}\right), 
\]
where $\xi$ is a parameter modeling the spatial correlation length.
In this case, if $\xi << 2 \pi$, the eigenvalues of the covariance operator satisfy
\[\lambda_k^2 = \exp\left(-\frac{\xi^2 k^2}{4 \pi}\right).
\] 
By substituting \eqref{5} into \eqref{6}, taking into account  \eqref{beta} and the ortogonality of the system $v_k$,
we obtain
\begin{equation} \label{7}
du_t^i=
\left[(I(x,t),v_i) - \alpha u_t^i + (KS)^i(\bar{u}_{t-\tau}) \right] dt + \epsilon \lambda_i d \beta_t^i,
\end{equation}
where  $(KS)^i(\bar{u}_{t})$ denotes the nonlinear term of the system:
\begin{equation} \label{8}
(KS)^i(\bar{u}_{t})= \int_{\Omega} v_i(x) \left( \int_{\Omega} K(|x-y|) S\left(\sum_{k=1}^{\infty} 
u_{t-\tau}^k v_k(y)\right)dy \right) dx. 
\end{equation}
Due to the symmetry of the kernel, the last integral may be rewritten as 
\begin{equation} \label{9}
(KS)^i(\bar{u}_{t-\tau})= \int_{\Omega} v_i (x) K(|x-y|) \left( \int_{\Omega}  S\left(\sum_{k=0}^{\infty} 
u_{t-\tau}^k v_k(y)\right)dy \right) dx. 
\end{equation}
When using the Galerkin method, we define an approximate solution by truncating the series expansion
\eqref{5}
\begin{equation} \label{10}
U_t^N(x)=\sum_{k=0}^{N-1} u_t^{k,N} v_k(x).
\end{equation}
Then the coefficients $u_t^{k,N}$ satisfy the following nonlinear system of stochastic delay differential
equations:
\begin{equation} \label{11}
du_t^{i,N}=
\left[(I(x,t),v_i) - \alpha u_t^{i,N} + (KS)^{i,N}(\bar{u}_{t-\tau})\right]dt  + \epsilon \lambda_i d \beta_t^i,
\end{equation}
where $(KS)^{i,N}(\bar{u}_{t})$ is given by 
\begin{equation} \label{8a}
(KS)^{i,N}(\bar{u}_{t-\tau})= h^2 \sum_{j=0}^{N}v_i(x) \left( \sum_{l=1}^{N} K(|x_l-x_j|) S\left(\sum_{k=1}^{N} u_{t-\tau}^k v_k(x_j)\right) \right) 
\end{equation}
$i=0,...,N-1$. In this case we are introducing in $[-L,L] $ a set of $N+1$ equidistant gridpoints 
$x_j=-l+j*h$, $j=0,...,N$, where $h=2l/N$, and using the rectangular rule to evaluate the integral
in \eqref{9}.

Since the problem has been reduced to  the system \eqref{11} (a system of nonlinear stochastic delay differential equations),
we can  apply the Euler-Maruyama method to the solution of this system. 
For the discretisation in time, we introduce on the interval $[0,T]$ a uniform mesh with step size $h_t$, such that $t_j= j h_t$, $j=0,1,...,n$. 
Then the solution $u_t^{k,N}$ of \eqref{11} will be approximated by a vector
 $(u_{1}^{k,N}, u_{2}^{k,N},..., u_{n}^{k,N})$, where
\[ u_{j}^{k,N} \approx u_{t_j}^{k,N} .\]
In these notations, the Euler-Maruyama method may be written as
\begin{equation} \label{13}
u_{j+1}^{i,N}=u_{j}^{i,N} + h_t \left[(I(x_i,t_j),v_i) - \alpha u_{j+1}^{i,N} + (KS)^{i,N}
(\bar{u}_{t_j-\tau} )\right] + \sqrt{h_t}\epsilon \lambda_i w_i,
\end{equation}
where $w_i$ is a random variable with normal distribution ($ w_i = N(0,1)$), $j=0,....,n$, $i=0,...,N-1$. 
 In the right-hand side of \eqref{13} we have written $\alpha u_{j+1}^{i,N}$, meaning that we are using a {\em semi-implicit} version of the Euler-Maruyama method.

Further we can rewrite equations \eqref{13} in the form
\begin{equation} \label{14}
u_{j+1}^{i,N}=\frac{ u_{j}^{i,N} + h_t \left[(I(x_i,t_j),v_i)+ (KS)^{i,N}(\bar{u}_{t_j} )\right] + \sqrt{h_t}\epsilon \lambda_i w_i}{1+\alpha h_t}.
\end{equation}
The meaning of the inner product in \eqref{14}  will be explained below.
In order to compute $u_{t_j-\tau}^k$ we should take into account that
$\tau= \frac{|x_{k_1}-x_{k_2}|} {v} $ (the time spent by the signal to  travel between  $x_{k_1}$
and $x_{k_2}$).
In general $\tau$ may not be a multiple of $h_t$.
Let $d$  and $\delta_t$ be the integer and the fractional part of  $\frac{\tau}{ h_t} $.
In this case, we have
\[ t_{j-d-1} \le  t_j -\tau  \le t_{j-d} 
\]
and
\[ h_t \delta_t= \tau-d h_t.
\]
The needed value of the solution  $u_{t_j-\tau}$  in \eqref{13} is then approximated by
\begin{equation} \label{inter}
  u_{t_j-\tau} \approx 
	\left\{ \begin{array} {c}
	u(t_{j-d}), \quad \makebox{if}  \,\delta_t <0.5, \\
	u(t_{j-d-1}), \quad \makebox{if} \,\delta_t \ge 0.5.
	\end{array}  \right.
\end{equation}
Concerning the choice of the basis functions, we consider a set of orthogonal functions, similar to the onde described in \cite{KR}, page 7.
 More precisely, we define
\begin{equation} \label{17}
v_k(x)=\exp (ikx), \qquad  k=0,1,..,N.
\end{equation}
Note that with this choice of the basis functions  the inner product in \eqref{14}  and the sums in
\eqref{8a} can
be interpreted as the Discrete Fourier Transform (DFT). In particular, the set of inner products
$(I(x,t_j),v_i) $ , $i=1,...,N$ may be seen as the DFT  of the vector $I_N$, which contains the values of the function $I(x,t)$ at the grid points   $x_k= -l + k h$, $ k=1,..,N$. In this case $t$ is fixed 
($t=t_j$).   Therefore, these inner products can be evaluated  efficiently by  the
 Fast Fourier Transform (FFT).

\section{Numerical Examples}
\subsection{Noise effects on stationary multi-bump solutions}
We have applied our algorithm to analyse the effect of noise on the formation of multi-bump solutions in dynamic neural fields, in the presence of space dependent external stimulli. In \cite{Laing}, the authors have investigated the formation
of regions of high activity (bumps) in neural fields, which can be switched 'on'
and 'off' by transient stimulli.
 The  stability analysis of such patterns in the deterministic case was carried out in
\cite{Flora} and in \cite{Erlhagen}. The effects of noise on such stationary pulse solutions was 
studied in \cite{kilpatrick}.  
 In this case,
the firing rate function $S(x)$ is the Heaviside function;
the connectivity kernel is given  by
\[ K(x)= 2\exp(-0.08 x ) \left( 0.08\sin(\pi x/10) +\cos(\pi x/10) \right).\]
In  \cite{Erlhagen} it was proved that such oscillatory connectivity kernels support the formation
of stationary stable multi-bump solutions, induced by external inputs. In this case the external input
has the form  
\[I(x) =-3.39967 +8 \exp\left( - \frac{x^2}{18}\right).\] 
This specific example was considered
in \cite{Flora}, p. 37. The parameters of the numerical approximation are $N=100, 
h=1, l=50; n=200, h_t=0.02.$
We start by considering the deterministic case. Using our code with $\epsilon=0$, we are able to reproduce   
three kinds of stationary solutions, which are displayed in Fig.1 and Fig.2. 
%%%%%%%%%%%%%  fig 1
\begin{figure}
\includegraphics[width=10cm]{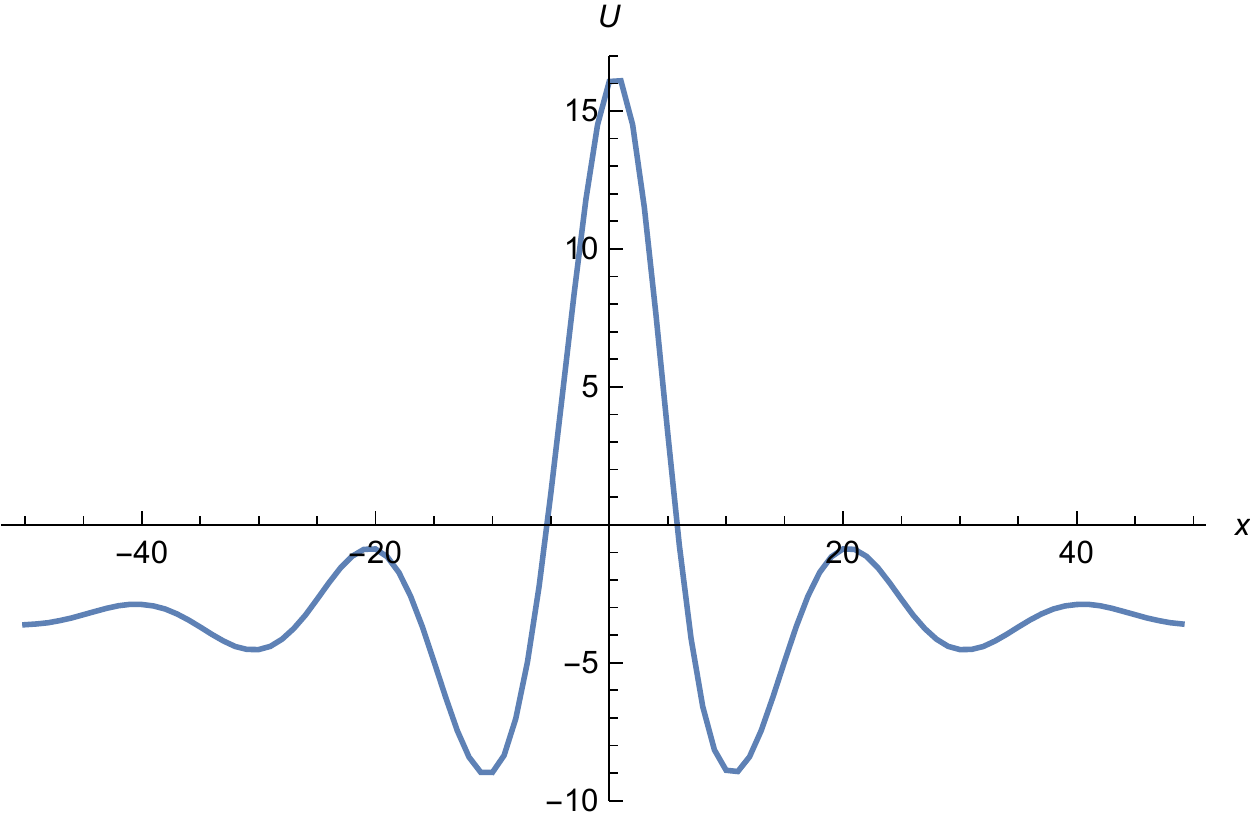} 
\caption{Deterministic case: stationary one-bump  solution}
\end{figure}
%%%%%%%%%%%%%  fig 2
\begin{figure}
$\begin{array}{cc}
\includegraphics[width =7cm]{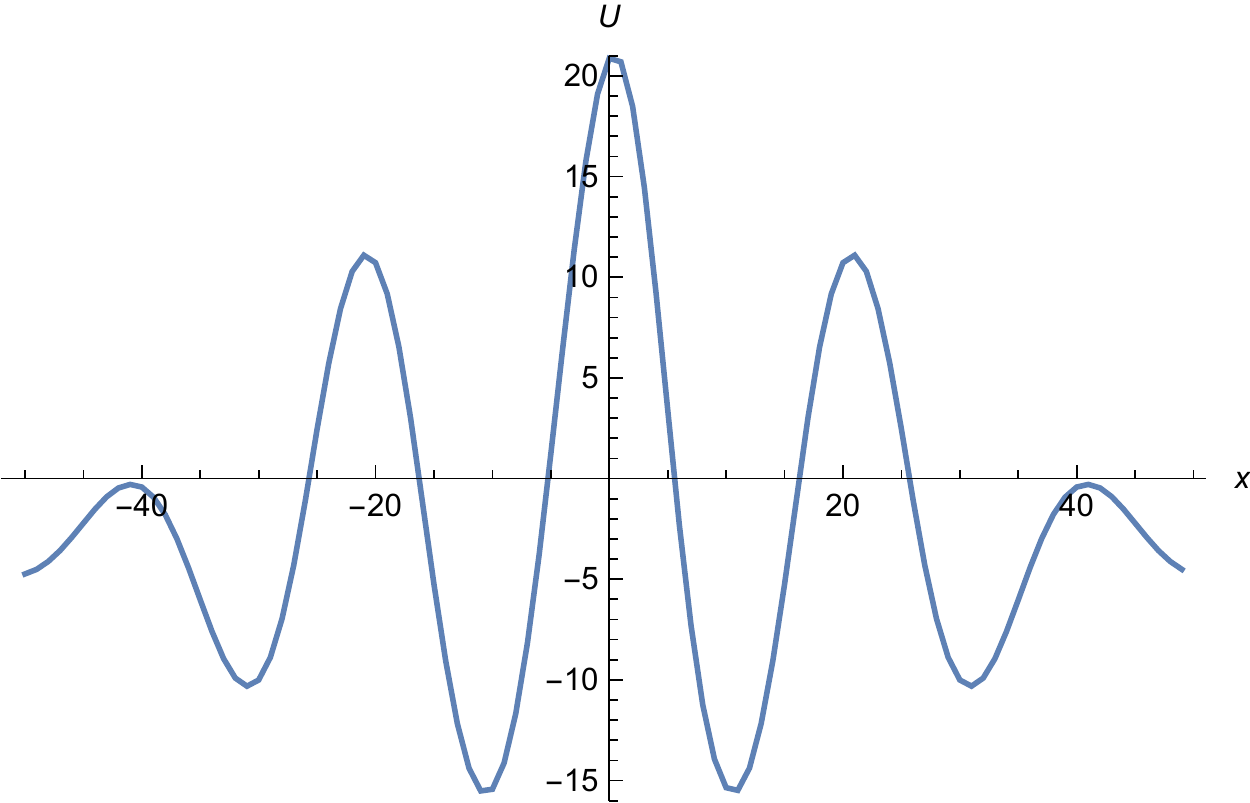} &
\includegraphics[width=7cm]{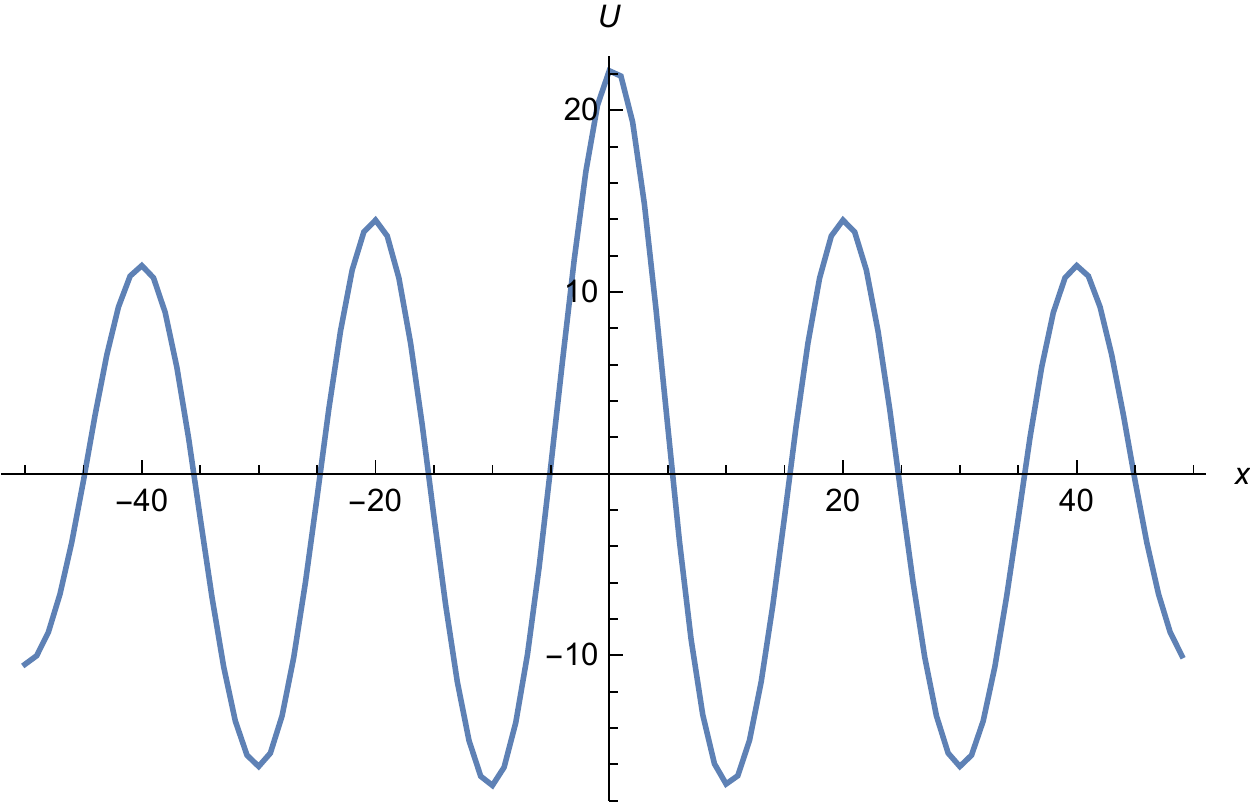} 
\end{array} $
\caption{ Deterministic case: three-bump (left) and five-bump (right) solutions}
\end{figure}
%%%%%%%%%%%%%%%%%%%%%%%%%%%%%%%%%%%%%%

As the {\bf first experiment} with the {\bf stochastic case}
we have performed a simulation with 100 paths, with  noise level
$\epsilon=0.01$ starting with the initial condition $U_0(x,t)\equiv 0$, over the time interval $t\in [0,4]$.
Our aim was to investigate the evolution of the paths of the stochastic equation and their relation with the stationary multi-bump solutions 
observed in the deterministic case.
We remark that each solution $V(x,t)$ of the deterministic equation (\ref{2}), as $t$ tends to infinity, tends to a certain
stationary solution $u(x)$, which is characterized by a certain maximal
value $u_{max} = \max_{x \in [-l,l]} u(x)$ (which is usually located at the origin) and a minimal value 
$u_{min} = \min_{x \in [-l,l]} u(x)$ (which is
attained at two symmetric points). Usually a solution may have other local minima and maxima, whose absolute values do not exceed 
those of $u_{min}$ and $u_{max}$.  The greater is the number of bumps, 
the higher is the value of $u_{max}$ and the lower is the value of $u_{min}$.
Note that $u_{max}$ and $u_{min}$ are respectively the limits, as $t \rightarrow \infty$,
of $u_{max}(t)$ and $u_{min}(t)$, that is, $u_{max}(t) = \max_{x \in [-l,l]} V(x,t)$
and $u_{min}(t) = \min_{x \in [-l,l]} V(x,t)$.
We have used these properties to analyse the paths of the stochastic equation.
 Let $u(s,x,t)$ denote the approximate value 
of $U_t(x)$, given by the $s$-th path. We use the following notations:
\[ U_{max,max}(t) = \max_{s\in\{1,\dots,100\}} \max_{i\in\{1,\dots,100\}} u(s,x_i,t);\] 
\[ U_{min,max}(t) = \min_{s\in\{1,\dots,100\}} \max_{i \in \{1,\dots,100\}} u(s,x_i,t);\] 
\[ U_{max,min}(t)= \max_{s\in\{1,100\}} \min_{i\in\{1,\dots,100\}} u(s,x_i,t);\] 
\[ U_{min,min}(t) = \min_{s\in\{1,\dots,100\} }\min_{i\in\{1,\dots,100\}} u(s,x_i,t);\] 
We consider the following approximations of mathematical expectations:
\[ E(u(x,t)) \approx \frac{1}{100} \sum_{s=1}^{100} u(s,x,t);\]
\[ E(\max_{x \in [-l,l]} u(x,t)) \approx E_{max}(t)= \frac{1}{100} \sum_{s=1}^{100} 
\max_{i \in \{1,\dots,100\}} u(s,x_i,t);\]
\[ E(\min_{x \in [-l,l]} u(x,t)) \approx E_{min}(t)= \frac{1}{100} \sum_{s=1}^{100} 
\max_{i \in\{1,\dots,100\}} u(s,x_i,t);\]
\begin{center}
%%%%%%%%%%%%%  fig 3
\begin{figure}
$\begin{array}{cc}
\includegraphics[width =7cm]{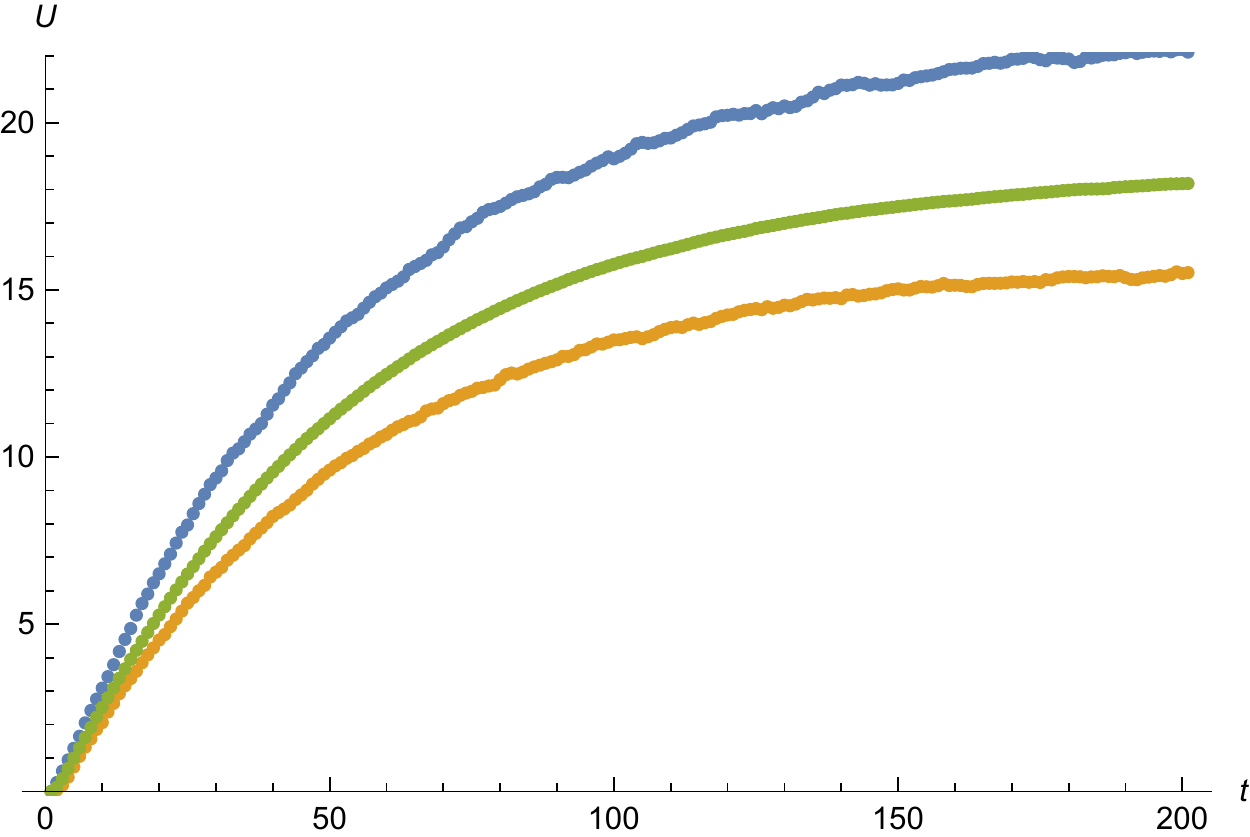} &
\includegraphics[width=7cm]{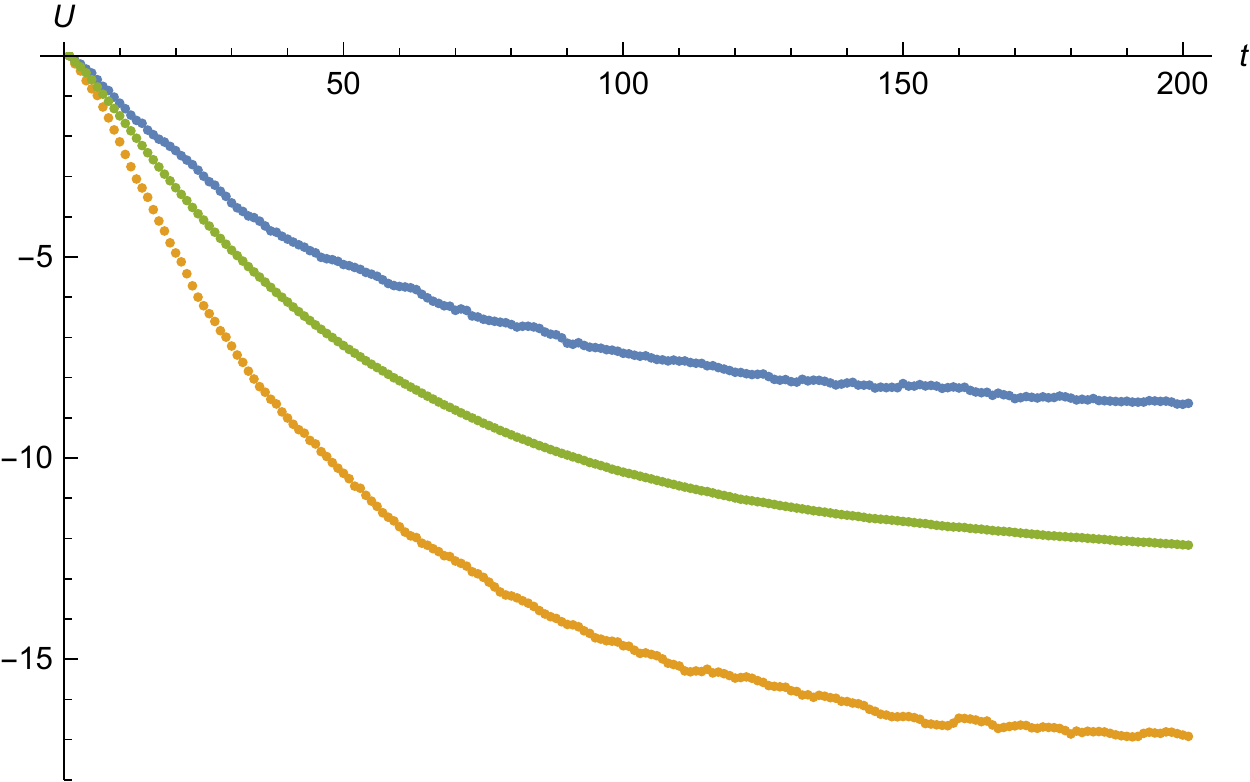} 
\end{array} $
\caption{Starting from the zero solution. Left: evolution of solution maximum - $U_{max,max}$ (blue), $U_{min,max}$
(yellow),  $E_{max}$ (green). Right: evolution of solution minimum:  $U_{max,min}$ (blue), $U_{min,min}$
(yellow),  $E_{min}$ (green). 
}
\end{figure}\end{center}
In Fig. 3  we can see the time evolution of the solution maximum (left) and the time evolution of the solution minimum (right). In the first case the graphs of  $U_{max,max}(t)$, $U_{min,max}(t)$ and $E_{max}(t)$ are displayed; the second graphic shows values of $U_{max,min}(t)$, $U_{min,min}(t)$ and
 $E_{min}(t)$ . These figures show that
the average value of the maximum increases with time, while the average value of the minimum decreases; as a result the average amplitude of the oscillations 
(height of the bumps) increases with time. Moreover, as it could be expected 
the dispersion (maximal difference between values of different paths) also increases with time. 
As it happens with $u_{min}(t)$ and $u_{max}(t)$ in the deterministic case, the values $E_{min}(t)$ and $E_{max}(t)$ stabilize as $t$ increases. 
For $j=200 (t_j=4)$, their values are close to the ones of $u_{min}$ and $u_{max}$, in the case of a deterministic one-bump solution.
%%%%%%%%%%%%%%%%%  fig 4
\begin{figure}
$\begin{array}{cc}
\includegraphics[width =7cm]{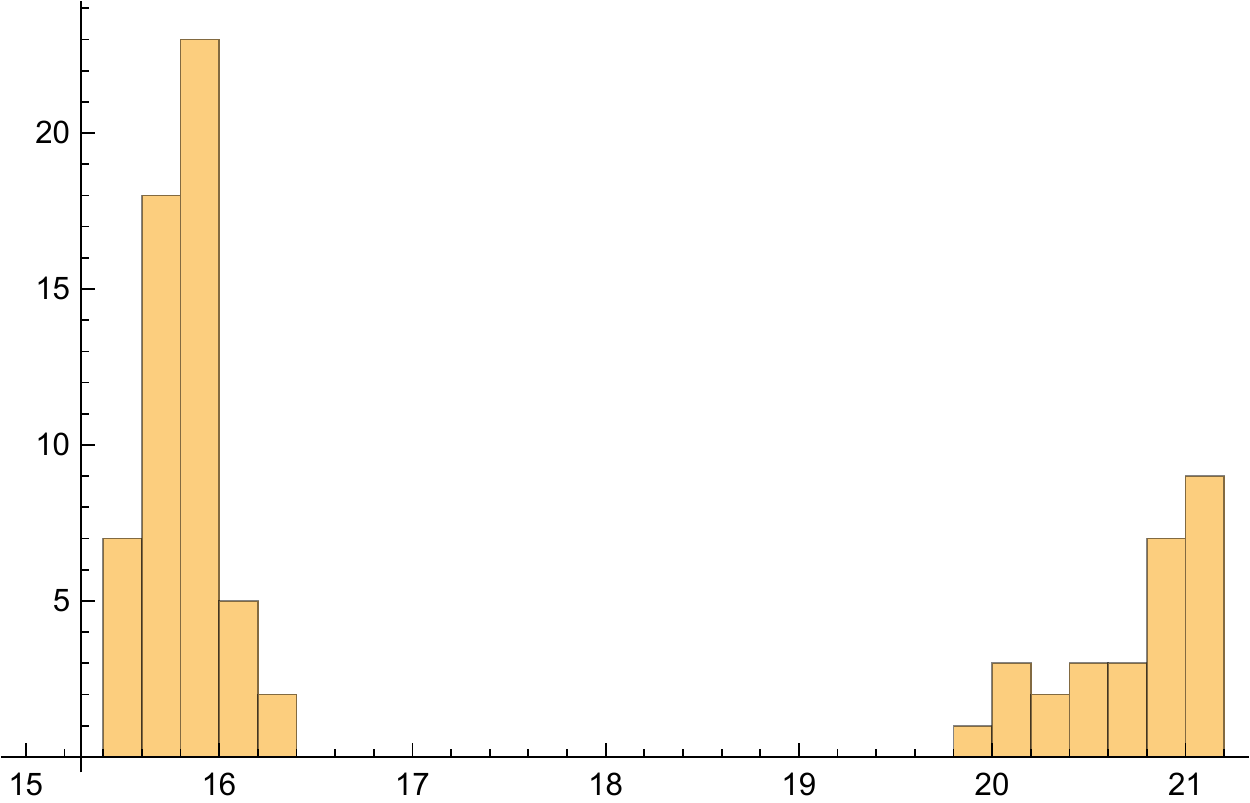} &
\includegraphics[width=7cm]{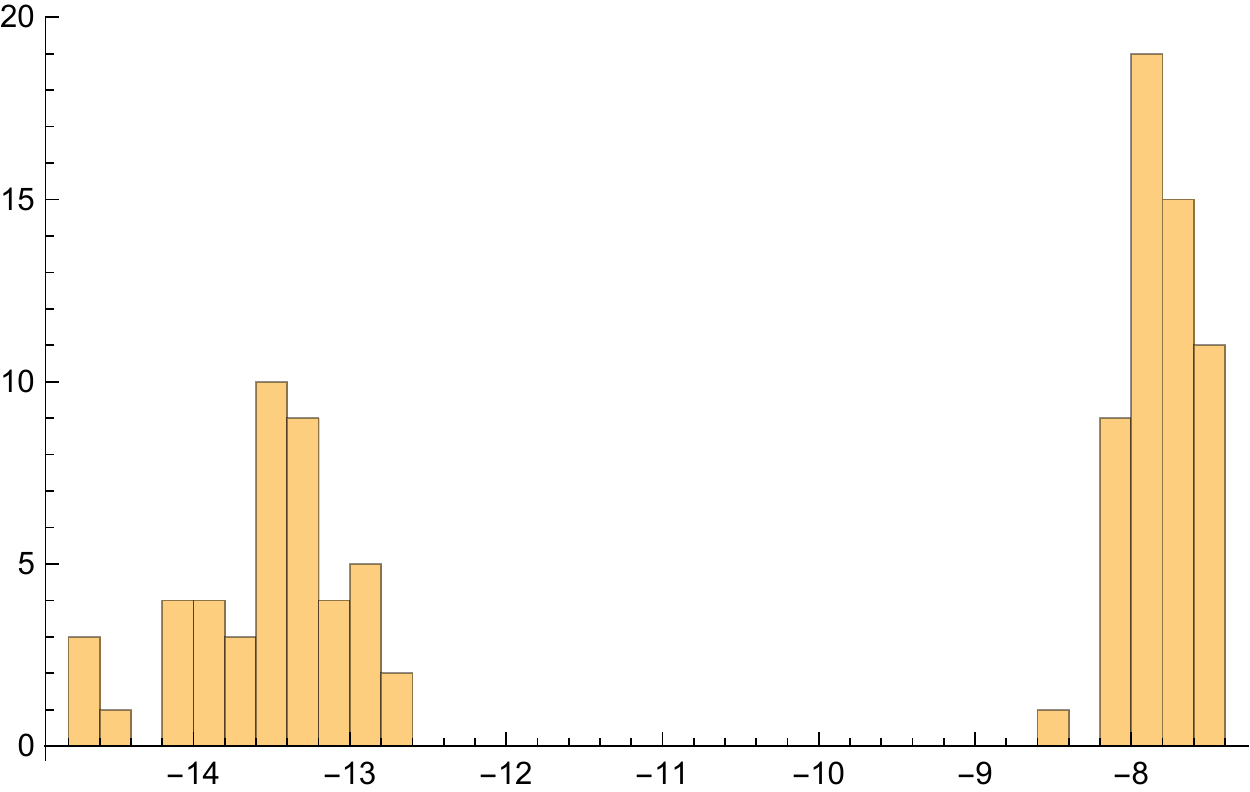} 
\end{array} $
\caption{ Histograms of distribution of $u_{max}$ (left) and 
$u_{min}$ (right), at $t_j=4$ ($j=200$), in the case $U_0(x,t)\equiv 0$, with $\epsilon=0.01$ .}
\end{figure}

Fig.4  shows the distribution (at $t=4$) of $u_{max}$ (left) and 
$u_{min}$ (right), at $t=4$,
for   the diferent paths. We see that the maxima are concentrated
 on the range $[15.8, 16.6]$, which corresponds to the stationary one-bump solution 
(see Fig. 1),  and 
$[20, 21.2]$,
which corresponds to three-bump and five-bump solutions (see Fig. 2). The minima
are concentrated  on the interval  $[-8.2,-7.4]$, which corresponds to the stationary one-bump solution,  and $[-14.,-12.5]$,
which corresponds to three-bump and five-bump solutions of the deterministic equation. This suggests that 
under the effect of not very strong noise, the most probable  values of the solution of the stochastic equation are close to the ones of the stationary solutions of the deterministic equation.

As a {\bf second numerical experiment}, we have performed a simulation with 100 paths, with  noise level
$\epsilon=0.01$, over the time interval $t\in [0,4]$, taking as initial condition $U_0(x,t)$ a  stationary one-bump solution (of the type represented in Fig.1).

%%%%%%%%%%%%%%%%%  fig 5
\begin{figure}
$\begin{array}{cc}
\includegraphics[width =7cm]{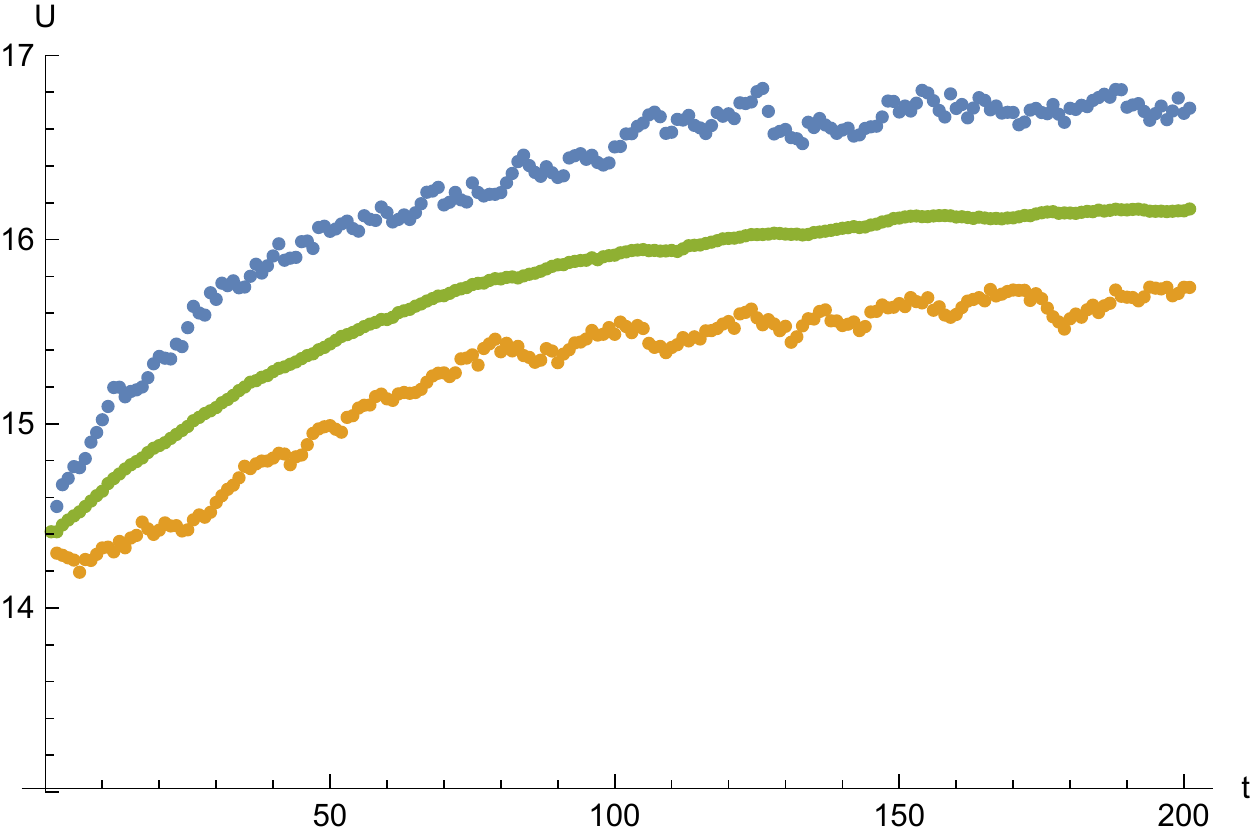} &
\includegraphics[width=7cm]{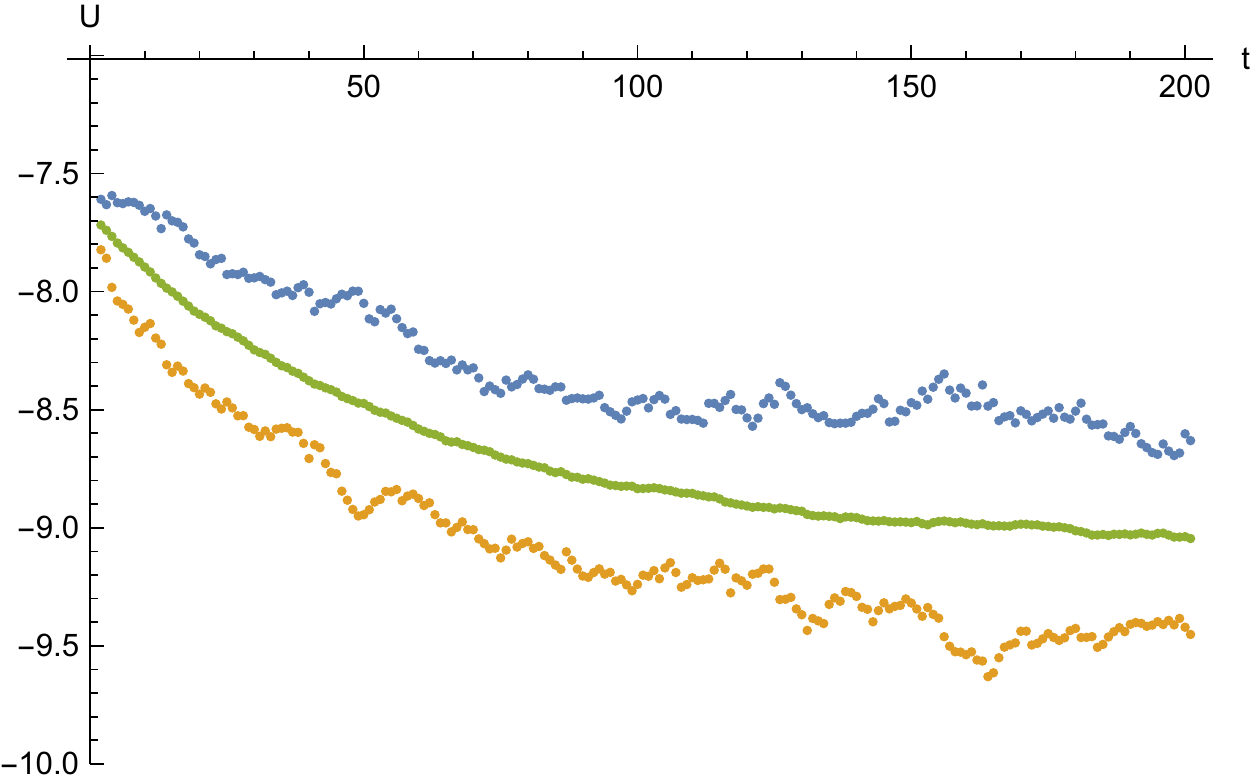} 
\end{array} $
\caption{Starting from a one-bump solution, with $\epsilon=0.01$. Left: evolution of solution maximum - $U_{max,max}$ (blue), $U_{min,max}$
(yellow),  $E_{max}$ (green). Right: evolution of solution minimum:  $U_{max,min}$ (blue), $U_{min,min}$
(yellow),  $E_{min}$ (green). 
}
\end{figure}
In Fig. 5 we again observe the evolution of $U_{max,max}(t)$, $U_{min,max}(t)$ and $E_{max}(t)$ (left-hand side); the values of values of $U_{max,min}(t)$, $U_{min,min}(t)$ and  $E_{min}(t)$ (right-hand side). 
 As in the previous case (Fig. 3) we observe that  
the average amplitude of the solutions oscillations increases with time; however for this level of noise ($\epsilon=0.01$), the amplitude of oscillations tends to stabilize and the mean value of the stochastic solution remains close to the one-bump deterministic solution.
%%%%%%%%%%%%%%%%%  fig 6
\begin{figure}
\includegraphics[width =10cm]{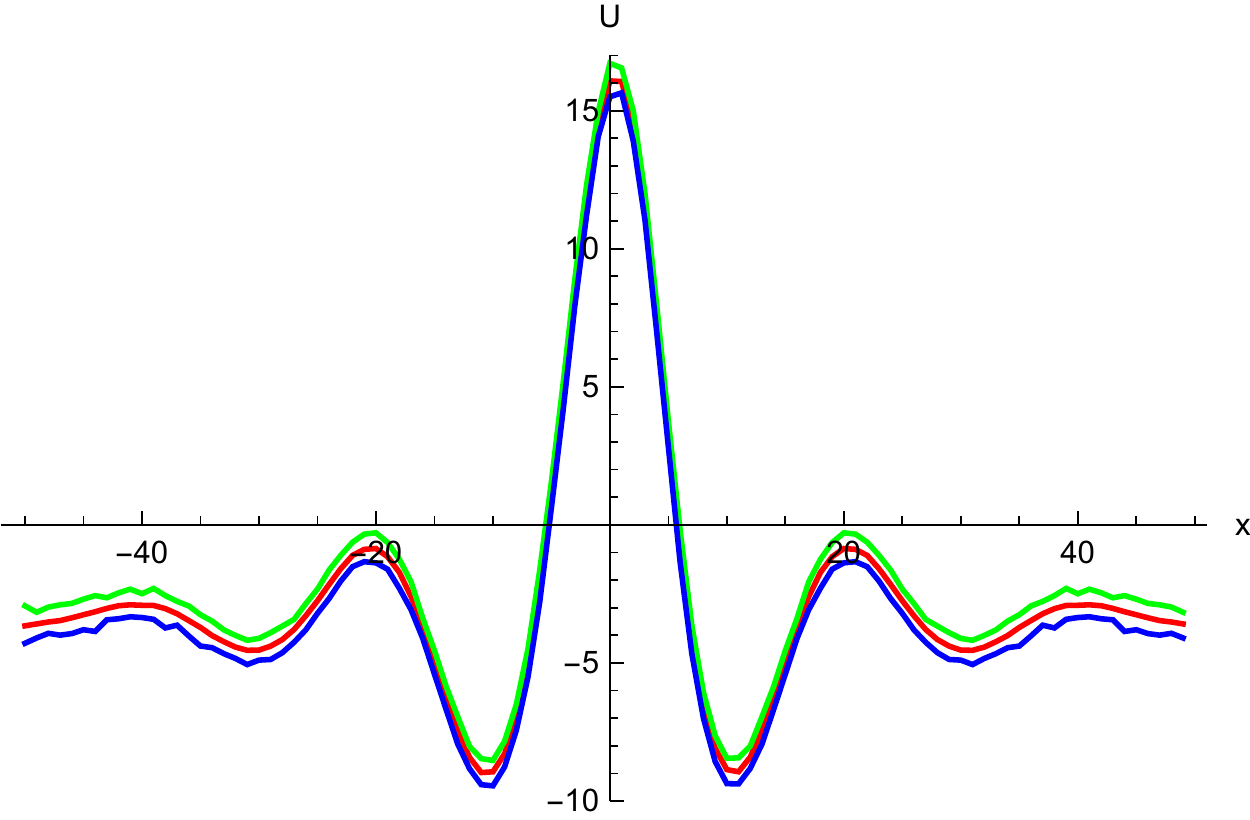} 
\caption{Graph of $ E(u(x,4))$ (red line),
$\max_{s \in\{1,\dots,100\}}u(s,x,4)$ (green line) and
$\min_ {s \in\{1,\dots,100\}} u(s,x,4)$ (blue line),} in the case $\epsilon=0.01$.
\end{figure}

In Fig. 6, the graph of the average solution $ E(u(x,t))$
at $t=4$  is
displayed, between the graphs containing the minimal and maximal values 
(of all the paths); these graphs correspond to the simulation starting with the stationary one-bump solution and with  noise level 
$\epsilon=0.01$.

%%%%%%%%%%%%%%%%%  fig 7
\begin{figure}
$\begin{array}{cc}
\includegraphics[width =7cm]{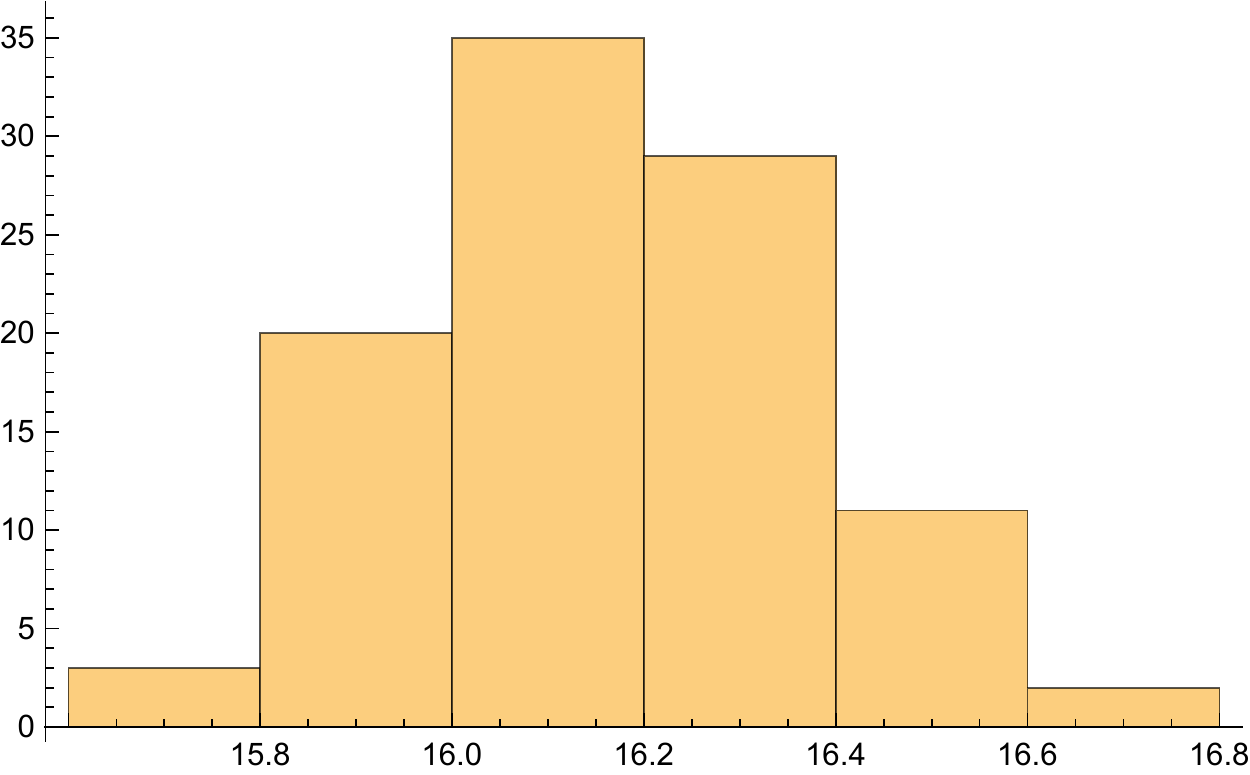} &
\includegraphics[width=7cm]{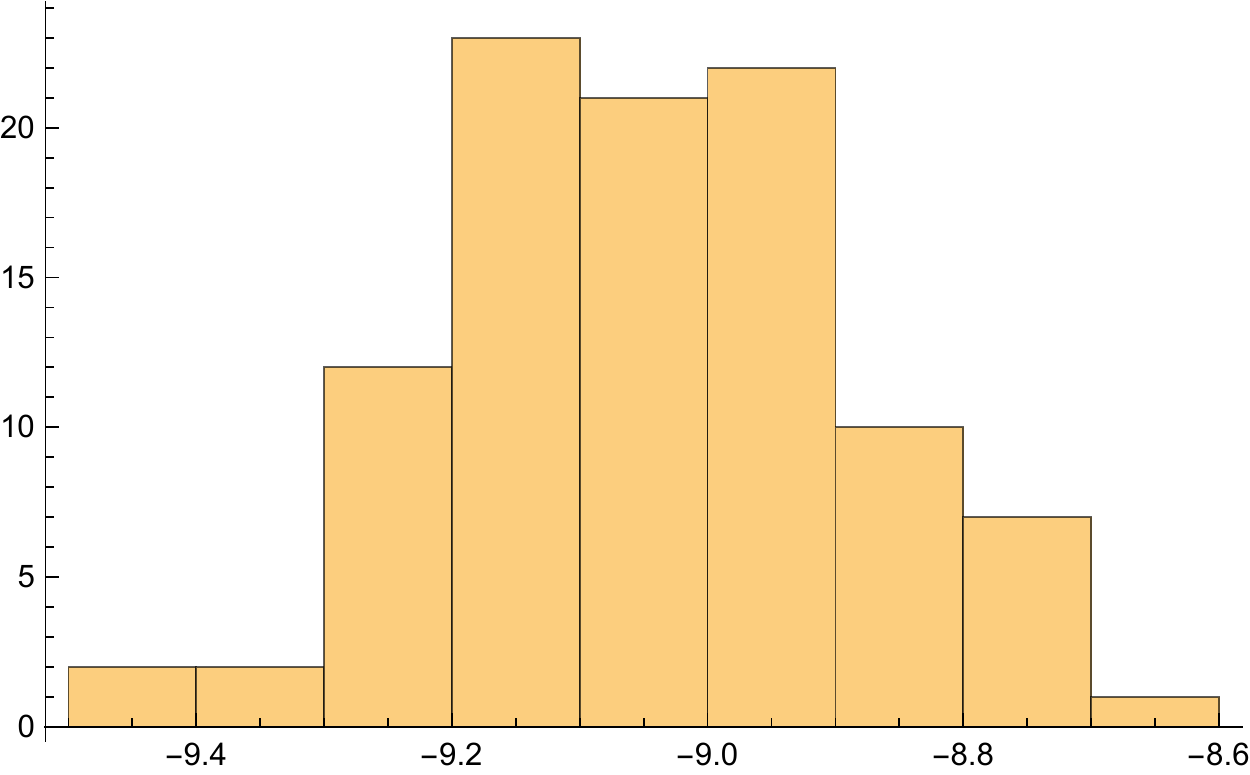} 
\end{array} $
\caption{Histograms of distribution of $u_{max}$ (left) and  $u_{min}$ (right), at $t=4$  (for the case where $U_0(x,t)$ is a deterministic one-bump solution, with $\epsilon=0.01$).}
\end{figure}
In Fig. 7 we can see the distribution (at $t=4$) of the values
of $u_{max}$ (left) and $u_{min}$ (right), for the same initial 
values and noise level.
  We see that the maxima are concentrated on the range $[15.8 ,16.6]$ and the minima are concentrated
on the range $[-9.4,-8.3]$. This means that
with level noise $\epsilon=0.01$ when the initial
value of the simulation is a deterministic one-bump solution,
the most probable values of the stochastic solution
are concentrated near the ones of  this stationary solution.

As a {\bf third numerical experiment}, we have performed a simulation with 100 paths, with level noise
$\epsilon=0.05$, over the time interval $t\in [0,4]$, taking again as initial condition $U_0(x,t)$ a  stationary one-bump solution.

%%%%%%%%%%%%%%%%%  fig 8
\begin{figure}
$\begin{array}{cc}
\includegraphics[width =7cm]{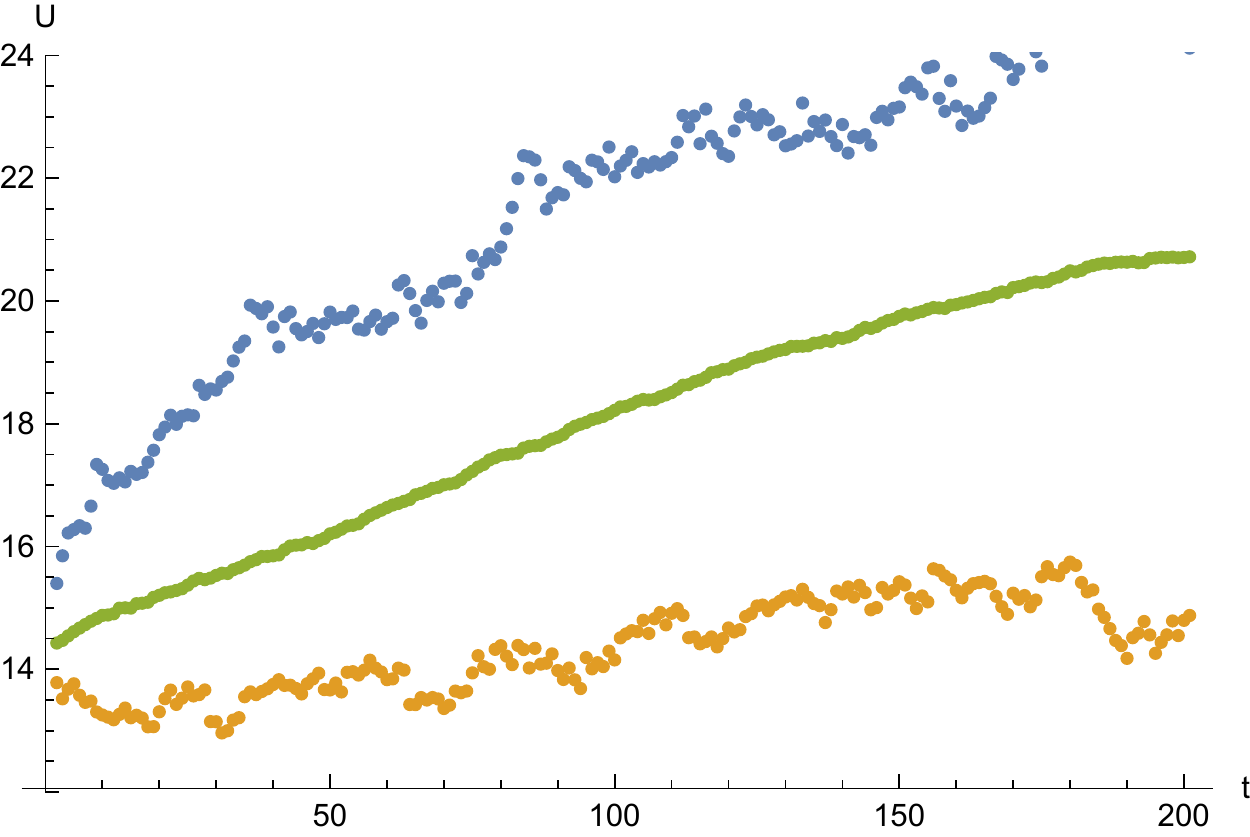} &
\includegraphics[width=7cm]{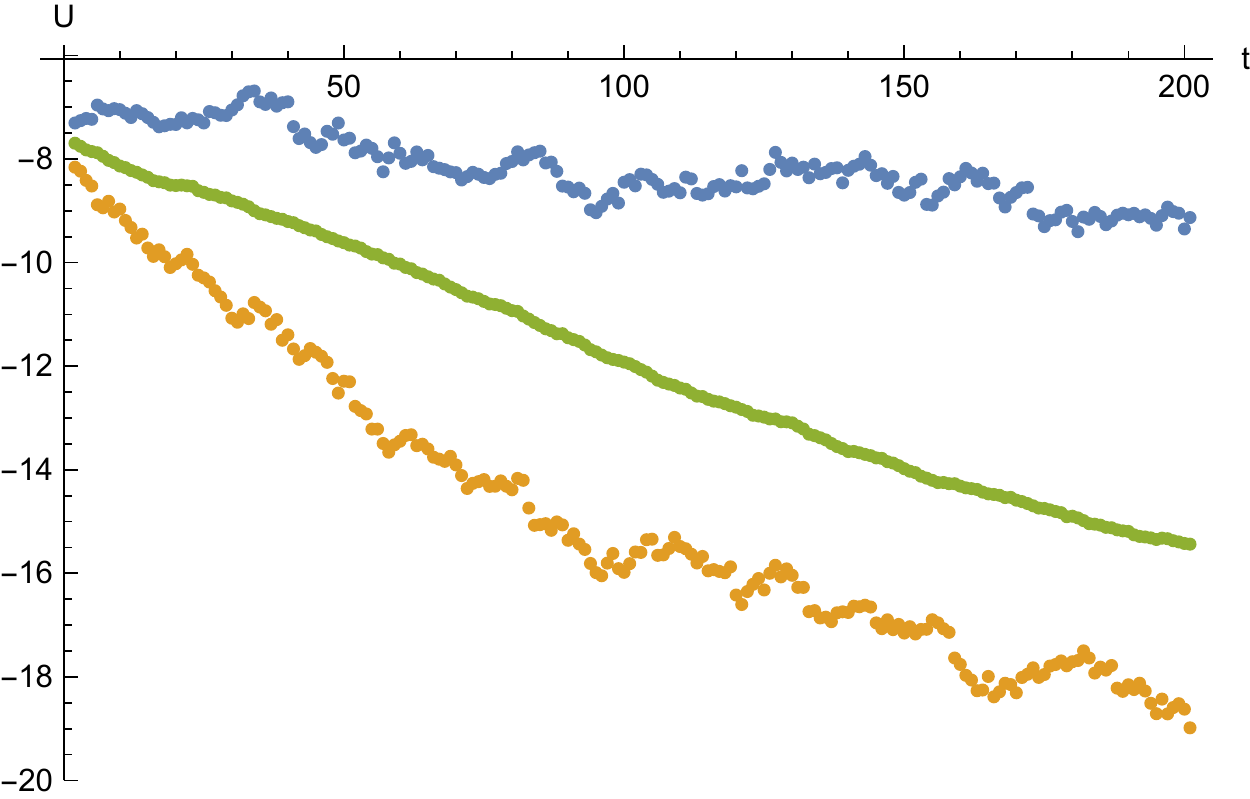} 
\end{array} $
\caption{Starting from a one-bump solution, with $\epsilon=0.05$. Left: evolution of solution maximum - $U_{max,max}$ (blue), $U_{min,max}$
(yellow),  $E_{max}$ (green). Right: evolution of solution minimum:  $U_{max,min}$ 
(blue), $U_{min,min}$ (yellow),  $E_{min}$ (green). }
\end{figure}
The evolution of the maximum and minimum
of the solutions is displayed in Fig. 8; in each case we can see the graphs of the average value of the 100 paths ($E_{max}(t)$ and  $E_{min}(t)$, respectively).  In Fig. 9, the graph of the average solution at $t=4$  is
plotted, between the graphs containing the minimal and maximal values (of all the paths). As in the previous cases we observe that  
the average amplitude of the solutions oscillations increases with time; moreover, for this level of noise ($\epsilon=0.05$), the amplitude of oscillations increases so much that some paths behave like three-bump
solutions  or five-bump solutions.

%%%%%%%%%%%%%%%%%  fig 9
\begin{figure}
\includegraphics[width =10cm]{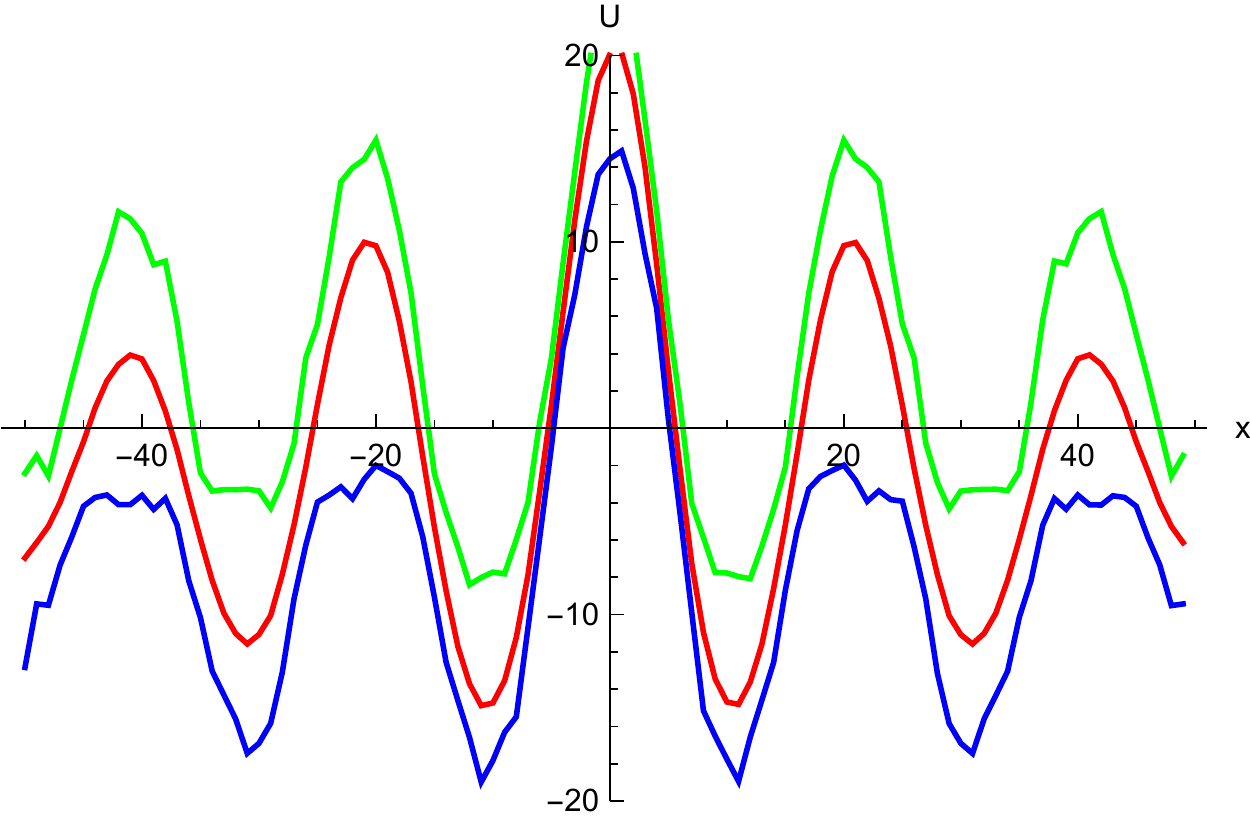} 
\caption{ Graph of $ E(u(x,4))$ (red line),
$\max_{s \in\{1,\dots,100\}}u(s,x,4)$ (green line) and
$\min_ {s \in\{1,\dots,100\}} u(s,x,4)$ (blue line), in the case $\epsilon=0.05$}.
\end{figure}

%%%%%%%%%%%%%%%%%%%%%%%%%%%%%%%%  fig. 10
\begin{figure}
$\begin{array}{cc}
\includegraphics[width =7cm]{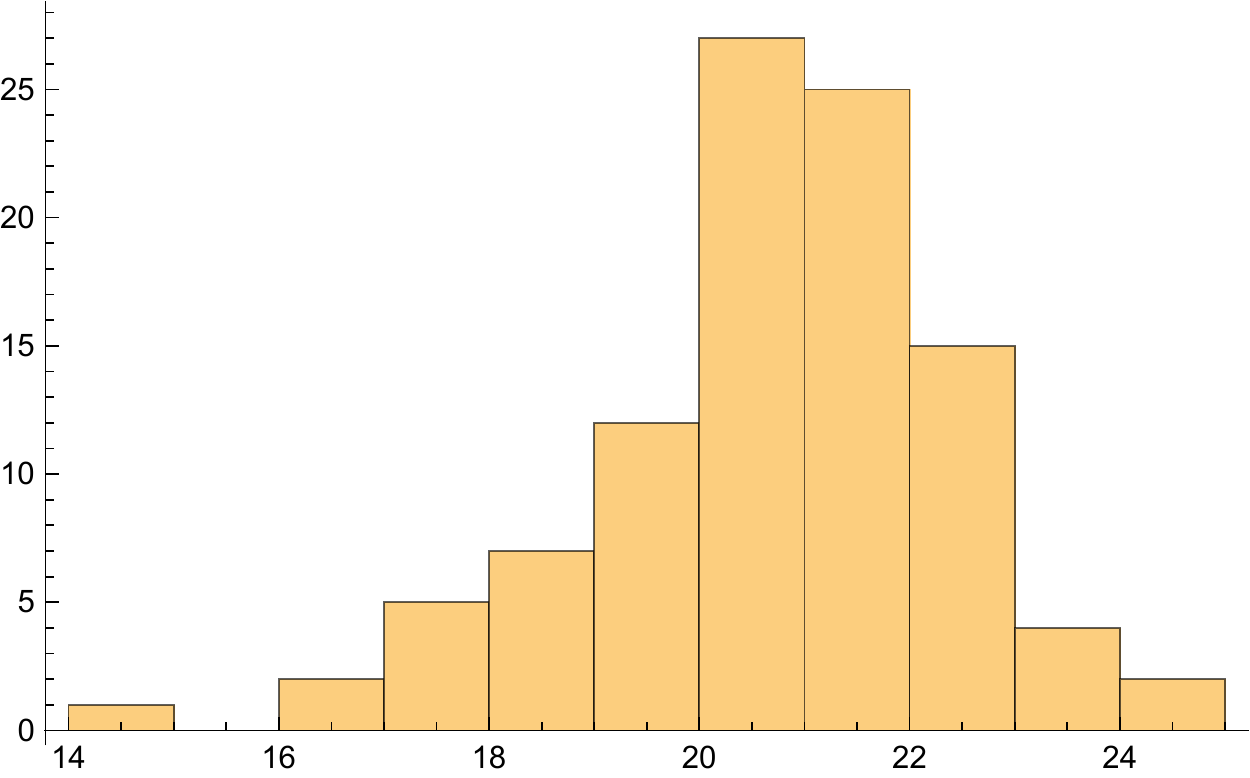} &
\includegraphics[width=7cm]{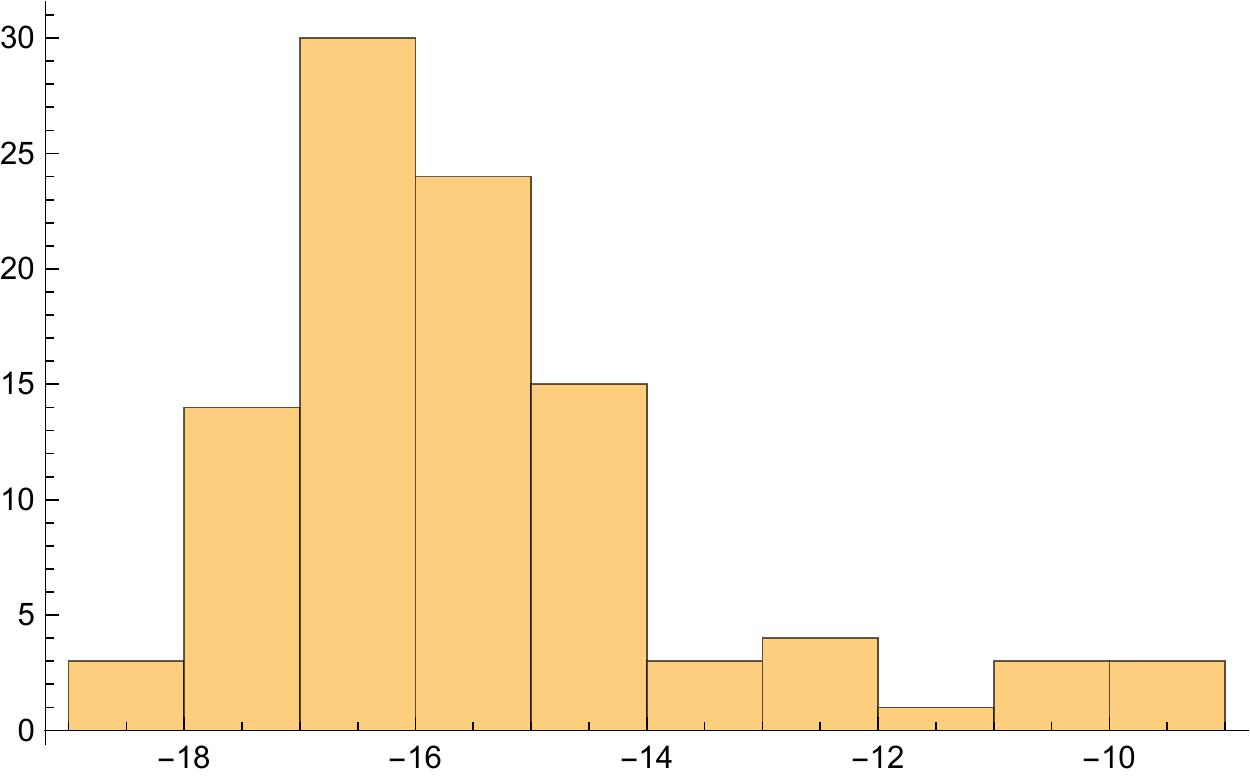} 
\end{array} $
\caption{Histograms of distribution of $u_{max}$ (left) and  $u_{min}$ (right), at
 $t=4$  (in the case where $U_0(x,t)$ is a deterministic one-bump solution, with
 $\epsilon=0.05$).}
\end{figure}

In Fig. 10  the distribution (at $t=4$) of the maxima and minima (respectively) of the diferent paths are displayed. We see that  most of the maxima are  on the range $[16,25]$ and  most of the minima are 
on the range $[-19,-9]$. This reflects the fact that
with  noise level $\epsilon=0.05$  many paths starting from a stationary one-bump solution can
after some time take values that are characteristic of three-bump or five-bump deterministic solutions.

The algorithm was implemented in Mathematica \cite{Wolfram} and 
the computations were performed in a PC with a 1.7Ghz processor
and 8 Gb of installed memory (RAM).
The computation of the numerical examples presented in this section (with 100 paths) takes about 1 hour.

\section{Conclusions}
In this paper we have introduced a new numerical algorithm for the approximation
of the stochastic neural field equation with delay. This numerical algorithm uses the
Galerkin method and  is inspired in the K\"{u}hn and Riedler's approach \cite{KR}.
The choice of the basis functions and grid points allows the use of the Fast Fourier
Transform to perform summations, which improves significantly the efficiency of the algorithm.

To test the algorithm we have applied it to the numerical solution of 
a neural field, in the presence of external stimulli, where stationary one-bump and multi-bump solutions are known to exist  in the deterministic case.
The numerical results suggest that for a low level of noise the trajectories of the stochastic equation
are concentrated near the stationary solutions of the deterministic one.  In particular, if the
initial condition is the trivial solution, in the deterministic case the solution tends with time
to a one-bump deterministic solution.
But in the presence of not very strong noise the trajectories of the stochastic equation split into
several classes, each of them close to  a different stationary solution of the deterministic equation.

In conclusion, we can say that our results are consistent with the study of Kilpatrick and Ermentrout \cite{kilpatrick},
where the authors conclude that upon breaking the translation symmetry of a neural field,
by introducing spatially heterogeneous input or synapses, bumps in the stochastic neural
field can become temporarily fixed to a finite number of locations.

\end{document}